\documentclass[12pt]{article}
\usepackage{amssymb,amsmath,amscd,amsthm}
\usepackage{graphicx,epsfig, color,float}

\usepackage{graphicx}
\usepackage[active]{srcltx}

\newtheorem{theorem}{Theorem}[section]

\newtheorem{lemma}[theorem]{Lemma}

\def\be{\color{black}}

\date{}

\begin{document}

\date{}
\title{Global limit theorem for parabolic equations with a potential}
\author{ L.
Koralov\footnote{Dept of Mathematics, University of Maryland,
College Park, MD 20742, koralov@math.umd.edu},
B. Vainberg\footnote{Dept of Mathematics and Statistics, UNCC, NC 28223,  brvainbe@uncc.edu}}
%\footnote{Dept of Mathematics, University of Maryland, College
%Park, MD 20742, koralov@math.umd.edu}
%\footnote{Dept of Mathematics, University of Maryland,
%College Park, MD 20742, mif@math.umd.edu}
\maketitle

\begin{abstract}  We obtain the asymptotics, as $t + |x| \rightarrow \infty$, of the fundamental solution to  the heat equation with a compactly supported potential. It is assumed that the corresponding stationary operator has at least one positive eigenvalue. Two regions with different types of behavior are distinguished: inside a certain conical surface in the $(t,x)$ space, the asymptotics is determined by the principal eigenvalue and the corresponding eigenfunction; outside of the conical surface, the main term of the asymptotics is a product of a bounded function and the fundamental solution of the unperturbed operator, with the contribution from the potential becoming negligible if $|x|/t \rightarrow \infty$. A formula for the global asymptotics,  as $t + |x| \rightarrow \infty$, of the solution in the entire half-space $t > 0$ is provided.

In probabilistic terms, the result describes the asymptotics of the density of particles in a branching diffusion with compactly supported branching and killing potentials.
\end{abstract}



\section{Introduction}

Let $p(t,x,y)$ be the fundamental solution of the parabolic equation with a continuous compactly supported potential $v$,
i.e., $p(\cdot, \cdot, y)$ satisfies
\begin{equation} \label{mpeq}
\frac{\partial}{\partial t} p(t,x,y) = \frac{1}{2} \Delta p(t,x,y) + v(x) p(t,x,y),~~~~t > 0,~ x \in \mathbb{R}^d,~~~~
p(0,x,y) = \delta_y(x).
\end{equation}
The spectrum of the operator $L  = \frac{1}{2} \Delta + v$ consists of the negative semi-axis $(-\infty, 0]$ (absolutely continuous spectrum) and at most a finite number of non-negative eigenvalues. We assume that it has at least one  positive eigenvalue. In this case, the largest eigenvalue $\lambda=\lambda_0 > 0$ is simple, and the corresponding eigenfunction (ground state) $\psi$ can be taken to be positive.
We also choose $\psi$ that is normalized: $\|\psi\|_{L^2(\mathbb{R}^d)} = 1$.

Our main result concerns the asymptotic behavior of $p(t,x,y)$ when $|y|$ is bounded and $t + |x| \rightarrow \infty$. We encountered the
problem when studying the distribution of particles in branching diffusions, but the result is of independent interest. In probabilistic terms, $p(t,x, y)$ is the density of particles at $x \in \mathbb{R}^d$ at time $t$ in a branching diffusion process that starts with a single particle located at $y \in  \mathbb{R}^d$. The particles undergo a Brownian motion and branching with
intensity $v$ in the regions where $v$ is positive; in the regions where $v$ is negative, the particles are killed at the rate $|v|$. The roles of $x$ and $y$ can be reversed:
$x$ can be viewed as the initial position of a particle, and $y$ can be the point where the density is observed.

It turns out that the interplay between the branching that takes place on the support of $v$ and the motion of the particles far away from the support of $v$ leads to the splitting of the space-time domain $[0,\infty) \times \mathbb{R}^d$ into two regions where $p$ has different types of asymptotic behavior. The regions are separated by the conical surface $C = \{(t,x): |x-y| = \sqrt{2 \lambda_0} t\}$ in $\mathbb R^{d+1}$ with the vertex at $x=y$.
We will formulate the main result (Theorem~\ref{mmt1}) in the interior and the exterior of the cone separately
(away from the boundary), and then (Part (c) of the theorem) will provide the asymptotics of $p$ valid in the entire region $t > 0$. The relation between the asymptotic formulas is discussed in the remarks following the theorem.

Let $\theta = \theta(t,x-y) = |x-y|/t$.  For $\varepsilon \geq 0$, let
\[
C^{\rm int}_\varepsilon = \{(t,x): \theta \leq  \sqrt{2 \lambda_0} -\varepsilon\},~~C^{\rm ext}_\varepsilon = \{(t,x): \theta \geq  \sqrt{2 \lambda_0} +\varepsilon\}.
\]

Throughout the paper, we will assume that $|y| \leq R$, where  $R > 0$ is a fixed constant such that ${\rm supp}(v)$ belongs to the ball $B_R$ of radius
$R$ centered at the origin. Let $p_0(t,x) = (2 \pi)^{-\frac{d}{2}} \exp(-{ |x|^2}/{2t})$ be the fundamental solution $p(t,x,0)$ corresponding
to $v \equiv 0$.   For $x \neq y$, let $\alpha = \alpha(x-y) = (x-y)/|x-y|$ be the unit vector in the direction of $x-y$. Recall that
\[
{\rm erf}(u) = \frac{2}{\sqrt{\pi}} \int_0^u e^{- s^2} ds,~~u \in \mathbb{R}.
\]
\begin{theorem} \label{mmt1}
(a) For each $\varepsilon > 0$, there is $\delta > 0$ such that
\begin{equation} \label{tt}
p(t,x,y) = e^{\lambda_0 t} \psi(x) \psi(y) (1  + O(e^{-\delta t}))~~~as~~t \rightarrow \infty,~(t,x) \in C^{\rm int}_\varepsilon, ~ 
|y| \leq R.  
\end{equation}
(b) For each $\varepsilon > 0$ and $(t,x) \in C^{\rm ext}_\varepsilon$,
\begin{equation} \label{tb}
p(t,x,y) =  p_0(t,x-y)
 (a(\theta, \alpha, y) + O(\frac{t}{|x-y|^2})~~~as~~|x| \rightarrow \infty,   ~|y| \leq R,  
\end{equation}
where $a$ is   a positive   continuous function equal to
\begin{equation} \label{taaa}
a(\theta, \alpha, y) = 1+ \int_0^\infty \int_{\mathbb{R}^d} e^{-\frac{\theta^2s}{2}-\theta\langle \alpha, y-z \rangle}v(z)
p(s,z, y)  d z ds,  ~|y| \leq R. 
\end{equation}
For large $\theta$, function $a$ has the following behavior:
\begin{equation} \label{asa}
a(\theta,\alpha,y) - 1 =\frac{1}{\theta}(\int_0^\infty v(y+s\alpha)ds+o(1)),  ~|y| \leq R,  
\end{equation}
where the remainder term can be replaced by $o(1/\sqrt\theta)$ if $v\in C^1$ or by $O(1/\theta)$ if $v\in C^2$.
(This implies that $a$ is bounded and separated from zero.)

(c) If $(t,x) \in (0,\infty) \times \mathbb{R}^d, ~  \theta \leq 1/\varepsilon_0$   with some
$\varepsilon_0>0$\be, and $t+ |x-y| \rightarrow \infty$,  then
\[
p(t,x,y)=e^{\lambda_0 t}\psi(x)\psi(y)(1+{\rm erf}(\sqrt{t}\frac{\sqrt { 2\lambda_0}-\theta}{\sqrt { 2 }}))(a_1(\theta,\dot{x},y)+O(|x|^{-1/2}))
\]
\begin{equation}\label{tc}
+p_0 (t,x-y)(a_2(\theta,\alpha,y)+O(|x-y|^{-1})),~~|y|\leq R,
\end{equation}
where $a_1,a_2$ are continuous (and, therefore, bounded) functions, and $a_1(\theta,\dot{x},y)={1}/{2}$ when $\theta\leq\sqrt{2\lambda_0}$.
\end{theorem}

\noindent
{\bf Remark 1.} For large $|x|$, the function $\psi$ can be replaced by its asymptotics at infinity:
\begin{equation} \label{ppsi}
\psi(x)=C(\dot{x})|x|^{\frac{1-d}{2}}e^{-\sqrt{2 \lambda_0}|x|}(1+O(|x|^{-1})), ~~~ |x|\to\infty,~~\dot{x} = x/|x|,
\end{equation}
with
\begin{equation} \label{intrep}
C(\dot{x}) = (2\pi)^\frac{1-d}{2}(2 \lambda_0)^\frac{d-3}{4} \int_{\mathbb{R}^d} e^{\sqrt{2 \lambda_0} \langle \dot{x}, z \rangle} v(z) \psi(z) dz   > 0. 
\end{equation}
The asymptotics (\ref{ppsi}) follows from the representation $\psi=(\Delta/2-\lambda_0)^{-1}(v\psi)$ and the explicit formula for the kernel $K$ of
the operator  $(\Delta/2-\lambda_0)^{-1}$.   Since
\[
\psi(x) = \int_{|z| < R} K(x-z) v(z) \psi(z) d z,
\]
and $K$ can be expressed through the Hankel function (see, e.g., formula (15) in \cite{CKMV}), one only needs to use the well-known asymptotics
of the Hankel function at infinity to get the integral representation (\ref{intrep}) for  $C$. The positivity of $C$ is a consequence of the maximum principle, which implies that $\psi\geq\psi_0,~|x|\geq R,$  where $\psi_0$ is the solution of the problem $(\frac{1}{2}\Delta-\lambda_0)\psi_0=0,~~|x|\geq R; ~~\psi_0=\min_{|x|=R}\psi$ on the boundary $|x|=R$.

Thus, formulas (\ref{tt}) and (\ref{ppsi}) imply that
\begin{equation} \label{pc}
p(t,x,y) = e^{\lambda_0 t -\sqrt{2 \lambda_0}|x| }   \psi(y)  C(\dot{x})|x|^{\frac{1-d}{2}} (1+O(|x|^{-1})) ~~{\rm as}~~ |x| \rightarrow \infty,~(t,x) \in C^{\rm int}_\varepsilon, ~
|y| \leq R.
\end{equation}
\\
{\bf Remark 2.}   Part (c) of the theorem provides the asymptotics of $p$ in the vicinity of $C$ that includes $C^{\rm int}_0$ and a part of $ C^{\rm ext}_0$.
% In fact, Part (c) gives the asymptotics in the entire
%half-space, as $|x| + t \rightarrow \infty$, which coincides with the asymptotics provided in Parts (a) and (b). \be
 The leading term of the asymptotics in Part  (c) coincides  with those in Parts (a) and (b) in $C^{\rm int}_\varepsilon$ and
$C^{\rm ext}_\varepsilon$, respectively.  This is discussed at the end of the proof and in Section~\ref{Att}.

It will be shown that,  in $C^{\rm int}_\varepsilon$, the second term in the right-hand side of (\ref{tc}) is
exponentially smaller (in $t$, as $t \rightarrow \infty$) than the first one.   The second term differs from the first one by a factor of order $O(1/\sqrt t)$ in a neighborhood of the cone $C$ where $\sqrt t|\theta  - \sqrt{2 \lambda_0}|\leq 1$. If  $0 < \theta - \sqrt{2 \lambda_0} \leq \gamma$ with a sufficiently small $\gamma$,  and $\sqrt t(\theta  - \sqrt{2 \lambda_0})\geq 1$, then the first term is of order $p_0(t,x-y)/(\theta  - \sqrt{2 \lambda_0})$, which implies that the first term is much larger than the second one if $\theta  - \sqrt{2 \lambda_0}$ is positive and small.
%the main terms of the asymptotics in (\ref{tc}) are equivalent when $\theta  - \sqrt{2 \lambda_0} $ is positive and separated from zero and %infinity, and the first term dominates when $\theta  - \sqrt{2 \lambda_0} $ is positive and small.
%For each $A > 0$, there is $\gamma > 0$ such that the first term
%is estimated from below by $A p_0(t,x-y)$ provided that $\theta  - \sqrt{2 \lambda_0} \leq \gamma$ and $t$ is sufficiently large.
If  $\theta  - \sqrt{2 \lambda_0} $ is positive and separated from zero and infinity, then both
terms are estimated from above, in absolute value, by $p_0(t,x-y)$. Their sum is estimated from below by a  positive multiple of $p_0(t,x-y)$ in such a region by Part (b) of the theorem. 

The proof of Part (c) is
partially based on Part~(a). Besides, Part (c) does not cover Parts (a) and (b) completely since the remainder terms in Parts (a) and (b) decay faster than the one in a neighborhood of the cone $C$. Part (b)  also provides a more precise description
of the remainder term in $C^{\rm ext}_\varepsilon$ and the asymptotics of the leading term as $\theta \rightarrow \infty$.  
% As $t + |x| \rightarrow \infty$, the function ${\rm erf}$
%in (\ref{tc}) is exponentially close to 1 in $C^{\rm int}_\varepsilon$ and to $-1$ in
%$C^{\rm ext}_\varepsilon$, and the leading terms in the asymptotic formulas (\ref{tt}) and (\ref{tb}) can be obtained from (\ref{tc}) by %using (\ref{ppsi}) and the asymptotics of the function ${\rm erf}$ for large values of the argument (see the end of the proof and %Section~\ref{Att}). The proof of Part (c) is
%partially based on Part~(a). Besides, Part (c) does not cover Part (b) completely since Part (b)  provides a more precise description
%of the remainder term in $C^{\rm ext}_\varepsilon$ and the asymptotics of the leading term as $\theta \rightarrow \infty$.
% From the asymptotics of ${\rm erf}$ and $\psi$, it follows that
\\
\\
{\bf Remark 3.} 
 As will be established in
the proof of the theorem, the integrand in (\ref{taaa}) decays exponentially in $s$, and so the value of the integral can be defined, with high accuracy,
by integration over a large but fixed interval in $s$.
\\

The difference in the types of asymptotic behavior inside and outside the cone is related to the following phenomenon, which we mention here just in passing. Let $N(t,x,U)$ be the
(random) number of particles in the branching diffusion process that are found in a domain $U \subseteq \mathbb{R}^d$ at time $t$, assuming that a single initial particle was located at $x \in \mathbb{R}^d$. Since $p$ can be interpreted as the density of particles,
\[
p(t,x,y) = \lim_{r \downarrow 0} \left({ \mathrm{E} N (t,x, B_r(y))}/{{\rm Vol} (B_r(y)}\right),
\]
where $B_r(y)$ is the ball of radius $r$ around $y$, and $\mathrm{E}$ stands for expectation.
For $x \in C^{\rm int}_\varepsilon$, the main contribution to $\mathrm{E} N(t,x, B_r(y))$ comes from the event that the initial particle
gets to the support of $v$. where branching occurs, much earlier than the given time $t$, undergoes branching, and the number of its descendants by time $t$ goes to infinity (the
probability of this event may go to zero, e.g., if $|x|$ grows linearly as a function of $t$).

For $x \in C^{\rm ext}_\varepsilon$, the main contribution to $\mathrm{E} N(t,x, B_r(y))$ comes from the event that the initial particle
reaches the support of $v$ at a time that is close to $t$, and the number of its descendants is bounded.
\\

Let us briefly discuss the relationship between our result and   asymptotic results for   branching diffusions
(see also \cite{GMV}, \cite{MV},  \cite{Bul}  for global asymptotics
for non-local operators and its applications to front propagation and intermittency). Consider a branching diffusion process that starts with a single particle. The questions concerning
front propagation and the structure of the population inside the front have been actively discussed in probabilistic and PDE literature.
The front can be defined as the boundary of the region $A_t$ occupied by particles, i.e., $x \in A_t$ if the probability of finding at
least one particle
at time $t$ in a unit neighborhood of $x$ exceeds  a fixed value $c \in (0,1)$. A somewhat different definition of the front is as the boundary of the region $B_t$, where $x \in B_t$ if the average number of particles at time $t$ in the unit neighborhood of $x$ exceeds a fixed positive constant.

As mentioned above, equation
(\ref{mpeq}) describes the density of particles in a branching diffusion, assuming that $v$ is the difference between the branching and
killing potentials. This immediately allows one to describe the evolution of $B_t$ using the asymptotics of the solution $p(t,x,y)$,
while the evolution of $A_t$ is described in terms of the solution to a related non-linear (FKPP) reaction-diffusion equation.

One of the first results on the front propagation is due to Bramson \cite{Bram}, who showed that
the front $\partial A_t$ lags by a logarithmic in $t$ distance behind $\partial B_t$ in the case of homogeneous branching.
This result has since been extended, including to the case of periodic branching, and refined (see, e.g., \cite{N1}, \cite{N2} and references therein). In the case of inhomogeneous branching, some of the foundational results (for the  leading term in speed of the FKPP front propagation) were
obtained by Freidlin (see, e.g.,   Chapter 6  of   \cite{Fr}). The analysis involves relating the solution of the FKPP equation to
the linear equation (\ref{mpeq}). In \cite{LS}, the logarithmic correction to the leading term for the front speed was studied in the case of a branching potential
(without killing) that is a sum of a constant and a rapidly decreasing function   (see also \cite{BH}, \cite{LS1}).   Whether the correction term appears or not depends on the relative strength of the constant (background) branching and the perturbation. In a forthcoming paper, we will use the global asymptotics of the solutions to linear equation  (\ref{mpeq}) and to similar equations on higher order correlation functions
of the particle field in order to provide a detailed  description of the particle field via the moment analysis. The structure
of the field, depending on the interplay between the localized branching/killing potential and the background potential, will be analyzed inside, at, and outside the front (in the latter case, the asymptotics of the probability to find a
particle near a given point is of interest). This approach based on the study of correlation functions was used in \cite{K1}, \cite{K2}  in the cases of constant or compactly supported branching potentials. Some of the recent results on the structure of the particle population in the case of rapidly decreasing branching potentials include \cite{S1}, \cite{S2},   \cite{NS}.  

Let us also mention that the results of the current paper can likely be generalized to case of periodic diffusion coefficients in equation
(\ref{mpeq}) (and  applied to  branching diffusions in periodic media with periodic branching/killing potential that is
perturbed by a compact function). The main difference is that here we use an explicit expression for the fundamental solution $p_0$ of
the unperturbed operator. In the periodic case, on the other hand, we have an asympotic formula (\cite{HK}), up to the pre-exponential term in the effective heat kernel, that is valid up to linear in time distances from the origin. Earlier results in this
direction are due to Norris \cite{No} and S. Agmon (unpublished). The asymptotics of Green's function for the corresponding elliptic
problem has also been studied extensively (see, e.g., \cite{MT},   \cite{KR}).  

\section{Proof of the main result}

\noindent
{\it Proof of Part (a).} A slightly different version of Part (a) was proved in \cite{K1}. We provide a simplified proof here.
First, we recall the simple arguments (see, e.g., \cite{CKMV},  Theorem 8.1) for the case when $|x| \leq R+1$.
For $\lambda \in \mathbb{C}' = \mathbb{C} \setminus (-\infty, 0]$, let $R_\lambda = (\frac{1}{2} \Delta + v - \lambda)^{-1}: L^2(\mathbb{R}^d)
\rightarrow L^2(\mathbb{R}^d)$ be the  resolvent of the operator $L  = \frac{1}{2} \Delta + v$. It is a meromorhpic function of $\lambda$ with
poles at eigenvalues of $L$. The norm of $R_\lambda$ does not exceed
the inverse distance of $\lambda$ from the spectrum, and therefore
\begin{equation} \label{dst}
\| R_\lambda \|_{L^2(\mathbb{R}^d)} \leq \frac{1}{|{\rm Im}(\lambda)|}.
\end{equation}

A similar estimate for large $|\lambda|$ is valid if $R_\lambda$ is viewed
as an operator from $C_0(K)$ to $C(\mathbb{R}^d)$ for a compact set $K \subset \mathbb{R}^d$ :
\begin{equation} \label{dst0}
\| R_\lambda f \|_{C(\mathbb{R}^d)} \leq \frac{c(K, \delta)}{| \lambda |}\|f\|_{C_0(K)}, \quad
|{\rm arg}\lambda|\leq\pi-\delta, ~~\delta>0,~~|\lambda| \rightarrow \infty.
\end{equation}
Indeed, for the unperturbed operator, estimate (\ref{dst0}) for the  resolvent $R_\lambda^0=(\Delta/2-\lambda)^{-1}$
follows (see \cite{CKMV}, Lemma 5.1) from the estimate on the kernel $R_\lambda^0(x-y)$ of operator $R_\lambda^0$:
\[
\int_K |R_\lambda^0(x-y)| dy \leq \frac{c(K)}{| \lambda |}, \quad  |{\rm arg}\lambda|\leq\pi-\delta, ~~\delta>0,~~|\lambda| \rightarrow \infty,
\]
which, in turn, is a simple consequence of the explicit formula for $R_\lambda^0(x-y)$. After that, (\ref{dst0}) follows from the resolvent identity: $R_\lambda=R_\lambda^0(I+v(x)R_\lambda^0)^{-1}$.

%Using the Laplace transform, we express the operator $P(t)$ with the integral kernel $p(t,x,y)$ as
%\begin{equation} \label{contour1}
%P(t) = \int_{{\rm Re}(\lambda) = \lambda_0 + 1} R_\lambda e^{\lambda t} d \lambda,
%\end{equation}
%and and
%\begin{equation} \label{contour1}
%p(t,x,y) = \int_{{\rm Re}(\lambda) = \lambda_0 + 1} R_\lambda(x,y) e^{\lambda t} d \lambda,
%\end{equation}
%where $R_\lambda(x,y)$ is the integral kernel of operator $R_\lambda$.
Let $\eta  = \eta(t,x)$ be a smooth function equal to zero when $t^2+|x|^2<1$ and equal to one when $t^2+|x|^2>2$. The function $p_\eta(t,x,y) = \eta(t,x-y) p(t,x,y)$ satisfies
\[
\frac{\partial p_\eta}{\partial t} = Lp_\eta+f,~~p_\eta|_{t = 0} = 0,
\]
where  
\[
f = p\frac{\partial \eta}{\partial t} - \langle   \nabla_x{\eta} ,\nabla_x{p} \rangle -
\frac{1}{2} p \Delta_x \eta.
\]
Observe  that $f=0$ when $t^2+|x-y|^2>2$ and that $f$ is infinitely differentiable for $t \geq 0$ and is equal zero together with all its
derivatives at $t = 0$.
 Using the Laplace transform, we obtain 
\[
\lambda \tilde{p}_\eta(\lambda,x,y)  = L \tilde{p}_\eta(\lambda,x,y) - \tilde{f}(\lambda,x,y),
\]
where $\widetilde{p}$  is the Laplace transform of $p$ and $\widetilde{f}$ is the Laplace transform of $f$.
Since $R_\lambda$ is analytic for ${\rm Re}(\lambda) > \lambda_0$,
\begin{equation} \label{contour1}
p_\eta(t,x,y) = -\int_{{\rm Re}(\lambda) = \lambda_0 + 1} R_\lambda \widetilde{f}(\lambda,\cdot,y) e^{\lambda t} d \lambda.
\end{equation}

Let $\varkappa$ be the distance from $\lambda_0 $ to the rest of the spectrum
of the operator $L$. The main term of the Laurent expansion of $-R_\lambda$ at the pole $\lambda_0$ is the operator with the integral kernel $\psi(x)\psi(y)/(\lambda - \lambda_0)$.
By (\ref{dst}), the contour of integration in (\ref{contour1}) can be shifted to the left, and therefore
\begin{equation} \label{contour2}
p_\eta(t,x,y) =e^{\lambda_0 t}\psi(x)\int_{\mathbb R^d}\psi(z)\widetilde{f}(\lambda_0,z,y)dz- \int_{{\rm Re}(\lambda) = \lambda_0 -\nu} R_\lambda \widetilde{f}(\lambda,\cdot,y) e^{\lambda t} d \lambda,
\end{equation}
where $\nu>0$ is an arbitrary positive number that is smaller than $\varkappa$.

Using integration by parts in the integral defining the Laplace transform of $f$, taking into account the properties of $f$ listed above, we obtain that for each $m > 0$,
 \begin{equation} \label{ef1}
 |\tilde{f}| \leq C_m |{\rm Im}(\lambda)|^{-m}~~~{\rm when}~ |x| \leq R+1, |y| \leq R.
 \end{equation}
  This estimate and (\ref{dst0}) imply that the
second integral above does not exceed $C e^{(\lambda_0 - \nu )t }$. Since $p(t,x,y) = p_\eta(t,x,y)$ for $t \geq 2$, it follows that
 \begin{equation} \label{inta}
p(t,x,y) =e^{\lambda_0 t}\psi(x)\int_{\mathbb R^d}\psi(z)\widetilde{f}(\lambda_0,z,y)dz+ O(e^{(\lambda_0 - \nu )t }),
 \end{equation}
when $|x| \leq R+1,~|y| \leq R,~t \rightarrow \infty$.

Let us show that the integral in (\ref{inta}) is equal to $\psi(y)$
%: \br
%\begin{equation} \label{ineql}
%\int_{\mathbb R^d}\psi(z)\widetilde{f}(\lambda_0,z,y)dz = \psi(y).
%\end{equation} \be
 Indeed,
\[
(L - \lambda) \widetilde{p}_\eta = \widetilde{f},~~~(L - \lambda) \widetilde{p} = \delta_y(x).
\]
Therefore, $\widetilde{f} = \delta_y(x) + (L - \lambda)(\widetilde{p}_\eta - \widetilde{p})$. It remains to substitute this expression with $\lambda = \lambda_0$ into the integral and note that the
term containing $(L - \lambda_0)$ vanishes since  $L$ is symmetric and  $(L-\lambda_0) \psi = 0$. We thus obtain
\begin{equation} \label{betaeq}
p(t,x,y) =e^{\lambda_0 t}\psi(x)\psi(y)+ \beta(t,x,y),
\end{equation}
where, for $\nu < \varkappa$,
\begin{equation} \label{betaeq1}
|\beta(t,x,y)| \leq  C(\nu)e^{(\lambda_0 - \nu )t },~~~|x| \leq R+1,~|y| \leq R,~t \rightarrow \infty.
\end{equation}

Let now $|x|>R+1$.  We demonstrated above that the first integral in the right-hand side
of (\ref{contour2}) is equal to $\psi(y)$. Using   the resolvent identity in the second integral,  we rewrite (\ref{contour2}) for $t \geq 2$ in the form
\begin{equation} \label{contour3}
p(t,x,y) =e^{\lambda_0 t}\psi(x)\psi(y)- \int_{{\rm Re}(\lambda) = \lambda_0 -\nu} R^0_\lambda (I+v(x)R^0_\lambda)^{-1} \widetilde{f}(\lambda,\cdot,y) e^{\lambda t} d \lambda,
\end{equation}
where $R^0_\lambda=(\frac{1}{2}\Delta-\lambda)^{-1}$ and the operator function
\[
(I+v(x)R^0_\lambda)^{-1}: ~L_2(\mathbb R^d)\to L_2(\mathbb R^d)
\]
is meromorphic in $C'$ with poles at eigenvalues of $L$. Since $v$ is compactly supported, this operator can also be viewed as operator in $L_2(B_R), ~B_R=\{x:|x|<R\}$. From (\ref{dst}) for $R^0_\lambda$, it follows that
\begin{equation} \label{nl2}
\|(I+v(x)R^0_\lambda)^{-1}\|_{L_2(B_R)} \leq {C}~~{\rm as}~|{\rm Im}(\lambda)| \rightarrow \infty.
\end{equation}
For the integral kernel of $R^0_\lambda$, the following estimate holds:
\[
|R^0_\lambda(x,z)| \leq C \frac{|e^{-\sqrt{2 \lambda}|x|}|}{|x|^{\frac{d-1}{2}}}|\lambda|^{\frac{d-3}{2}},~~~|z| \leq R,~|x| \geq R +1.
\]
Thus, for $\lambda$ such that ${\rm Re}(\lambda) = \lambda_0 - \nu$,
\[
|R^0_\lambda(x,z)| \leq C \frac{e^{-\sqrt{2 (\lambda_0-\nu)}|x|}}{|x|^{\frac{d-1}{2}}}|\lambda|^{\frac{d-3}{2}},~~~|z| \leq R,~|x| \geq R +1.
\]
From here, (\ref{ef1}) and (\ref{nl2}) it follows that the following estimate is valid for the integral $I(\nu)$ in (\ref{contour3})
\[
|I(\nu)|\leq C|x|^{\frac{1-d}{2}}e^{(\lambda_0-\nu)t-\sqrt{2 (\lambda_0-\nu)}|x|}=C|x|^{\frac{1-d}{2}}e^{\lambda_0t-\sqrt{2 \lambda_0}|x|+t\zeta},
\]
where $\zeta=-\nu+(\sqrt{2 \lambda_0}-\sqrt{2 (\lambda_0-\nu)})|x|/t$. Since $\theta=|x-y|/t\leq \sqrt{2 \lambda_0}-\varepsilon$ in $C^{\rm int}_\varepsilon$, it follows that $|x|/t\leq \sqrt{2 \lambda_0}-\varepsilon/2$ when $(t,x)\in C^{\rm int}_\varepsilon$ and $t$ is large enough. Thus, for those values of $(t,x)$, $\zeta$ does not exceed $-\delta$ with
\[
\delta=\nu-(\sqrt{2 \lambda_0}-\sqrt{2 (\lambda_0-\nu)})(\sqrt{2 \lambda_0}-\frac{\varepsilon}{2})=\nu(1-\frac{2(\sqrt{2 \lambda_0}-\varepsilon/2)}{\sqrt{2 \lambda_0}+\sqrt{2 (\lambda_0-\nu)}}),~~\nu < \varkappa. 
\]
For each $\varepsilon>0$ we can choose $\nu=\nu(\varepsilon)>0$ so small that $\delta=\delta(\varepsilon)>0$. Then, for sufficiently large $t$,
\[
|I(\nu)|\leq C|x|^{\frac{1-d}{2}}e^{(\lambda_0-\delta)t-\sqrt{2 \lambda_0}|x|}~~~{\rm in} ~~C^{\rm int}_\varepsilon.
\]
From (\ref{ppsi}) and positivity of $C(\dot x)$, it follows that
the right-hand side in the estimate on $I(\nu)$ can be replaced by $C\psi (x)e^{(\lambda_0-\delta)t}$. This, together with (\ref{contour3}), completes the proof of Part (a) for $|x| > R+1$. It remains to deal with Parts (b) and (c).

The analysis will be based on the Duhamel formula:
\begin{equation} \label{duh}
p(t,x,y) = p_0(t,x-y) +I,
\end{equation}
\begin{equation} \label{II}
I:= \int_{\mathbb{R}^d}\int_0^t\frac{1}{(2 \pi (t-s))^{d/2}} e^{-\frac{|x-z|^2}{2(t-s)}} v(z) p(s,z,y)ds dz.
\end{equation}
\\

\noindent
{\it Proof of Part (b).}
We will need the following two relations. Since $p(t,x,y)<C(T)p_0(t,x,y)$ on any bounded time interval $0<t\leq T$ and $p_0>Ct^{-d/2}$
when $t\geq 1,~|x|\leq R$, from Part (a) of the theorem it follows that, for arbitrary $\widetilde{\lambda}_0 > \lambda_0 $ and some $C=C(\widetilde{\lambda}_0)$, the following estimate holds
\begin{equation} \label{estp}
 |v(z)|p(s,z,y) \leq C |v(z)|p_0(s,z-y)e^{\widetilde{\lambda}_0 s},  ~~s \geq 0,~z \in \mathbb{R}^d, ~|y| \leq R. 
\end{equation}
The second relation concerns an expansion of the exponential factor in the integrand of $I$ that leads to a
representation of $I$. We have
\[
\frac{|x-z|^2}{2(t-s)}=\frac{|x-y|^2}{2(t-s)}+\frac{\langle x-y,y-z\rangle}{t-s}+\frac{|y-z|^2}{2(t-s)}.
\]
We represent $1/(t-s)$ in the first term as $1/t+s/t^2+s^2/[t^2(t-s)]$ and in the second term as $1/t+s/[t(t-s)]$. This leads to
\[
\frac{|x-z|^2}{2(t-s)}=
\frac{|x-y|^2}{2t}+\frac{\theta^2s}{2}+\frac{\theta^2s^2}{2(t-s)}+\theta\langle \alpha, y-z \rangle+\frac{\theta\langle
\alpha, y-z \rangle s}{t-s}+\frac{|y-z|^2}{2(t-s)}
\]
\[
=\frac{|x-y|^2}{2t}+\frac{\theta^2s}{2}+\theta\langle \alpha, y-z \rangle+\frac{|\theta s\alpha+y-z|^2}{2(t-s)}.
\]
Hence, $I$ can be rewritten in the form
\begin{equation}\label{ib}
I=p_0(t,x-y)\int_{\mathbb{R}^d}\int_0^t(\frac{t}{ t-s})^{d/2} e^{-\frac{\theta^2s}{2}-\theta\langle \alpha, y-z \rangle-\frac{|\theta s\alpha+y-z|^2}{2(t-s)}} v(z) p(s,z,y)ds dz.
\end{equation}
The remainder term in (\ref{tb}) (let us denote it by $r$) equals
\[
r=\int_{\mathbb{R}^d}\int_0^t [(\frac{t}{ t-s})^{d/2}e^{-\frac{|\theta s\alpha+y-z|^2}{2(t-s)}}-1]e^{-\frac{\theta^2s}{2}-\theta\langle \alpha, y-z \rangle}v(z) p(s,z,y)ds dz
\]
\[
-\int_{\mathbb{R}^d}\int_t^\infty e^{-\frac{\theta^2s}{2}-\theta\langle \alpha, y-z \rangle}v(z) p(s,z,y)ds dz =: r_1 - r_2.
\]

Let us estimate the remainder. We will start with $r_1$. Using (\ref{estp}) and the relation
\begin{equation}\label{rr2}
-\frac{\theta^2s}{2}-\theta\langle \alpha, y-z \rangle-\frac{|y-z|^2}{2s}=-\frac{|\theta s\alpha+y-z|^2}{2s},
\end{equation}
followed by substitution $z=y+\theta s\alpha+\sqrt su$, we  obtain
\[
|r_1|\leq C\int_{\mathbb{R}^d}\int_0^t |(\frac{t}{ t-s})^{d/2}e^{-\frac{|\theta s\alpha+y-z|^2}{2(t-s)}}-1|\frac{1}{s^{d/2}}e^{-\frac{|\theta s\alpha+y-z|^2}{2s}+\widetilde{\lambda}_0s}|v(z)|ds dz
\]
\begin{equation}\label{rr}
=C\int_{\mathbb{R}^d}\int_0^t |(\frac{t}{ s(t-s)})^{d/2}e^{-\frac{|\theta s\alpha+y-z|^2t}{2s(t-s)}}-\frac{1}{s^{d/2}}e^{-\frac{|\theta s\alpha+y-z|^2}{2s}}|e^{\widetilde{\lambda}_0s}|v(z)|ds dz
\end{equation}
\[
=C\int_{\mathbb{R}^d}\int_0^t |(\frac{t}{ t-s})^{d/2}e^{-\frac{|u|^2t}{2(t-s)}}-e^{-\frac{|u|^2}{2}}|e^{\widetilde{\lambda}_0s}|v(y+\theta s\alpha+\sqrt su)|ds du.
\]

Function $v$ above vanishes if $|\theta s\alpha|>|\sqrt su|+|y|+R$, i.e., the integration in $s$ can be performed over those values for which $\theta s<\sqrt s |u|+2R$. By solving this quadratic inequality for $\sqrt s$, we obtain that $\sqrt s<(|u|+\sqrt {|u|^2+8R\theta})/2\theta$, and therefore,
\begin{equation}\label{4}
s<s_0:=\frac{2(|u|^2+(\sqrt {|u|^2+8R\theta})^2)}{4\theta^2}=\frac{|u|^2+  4R\theta}{\theta^2}.
\end{equation}

Note that, for $(t,x) \in C^{\rm ext}_\varepsilon$ and  $|u|$ bounded by a fixed constant, we have
\begin{equation}\label{s0t}
\frac{s_0}{t}=\frac{|u|^2t}{|x-y|^2}+\frac{ 4R }{|x-y|}<  \frac {C(1+|u|^2)}{|x-y|}\ll 1  \quad {\rm  when} ~~t^2 +| x-y|^2\to\infty.
\end{equation}
This allows us to estimate the difference of exponents (with the pre-exponential factor) in the
integrand in the right-hand side of (\ref{rr}). We denote this difference by $w$:
\[
w=w(b)=b^{-d/2}e^{-\frac{|u|^2}{2b}}-e^{-\frac{|u|^2}{2}}, \quad b:=1-\frac{s}{t}.
\]
Then,  by the mean value theorem,  $w=-w'(b^*)s/t$, where 
\[
w'(b) = (\frac{u^2}{2 b^2} - \frac{d}{2b}) b^{-d/2}e^{-\frac{|u|^2}{2b}}
\] 
is the derivative of $w$ in $b$ and $b^*\in[b ,1]$.
If $|u|>\sqrt{d}$, then $w'(b)>0$ for $0\leq b\leq1$, and therefore $|w'(b^*)|\leq w'(1)$. If $|u|\leq \sqrt{d}$, then $s_0/t \ll 1$, and therefore $b \in [1/2,1]$. This implies that
$|w'(b^*)|<C(1+|u|^2)e^{-|u|^2/2}$. Hence, the latter estimate is valid for all $|u|$. From here, (\ref{rr}), (\ref {s0t}) and the boundedness of $|v| $ it follows (using integration in $s$) that
\[
|r_1|\leq \frac {C}{|x-y|}\int_{\mathbb{R}^d}\int_0^{s_0} (1+|u|^2)^2e^{-|u|^2/2}e^{\widetilde{\lambda}_0 s}ds du
\]
\[
= \frac {C}{\widetilde{\lambda}_0 |x-y|}\int_{\mathbb{R}^d}(1+|u|^2)^2  e^{-|u|^2/2}  [ e^{\widetilde{\lambda}_0\frac{|u|^2+4R\theta}{\theta^2}} -1 ] du.
\]
Since $\theta>\sqrt{2\lambda_0}$ in ${C^{\rm ext}_\varepsilon}$, we can choose
$\widetilde{\lambda}_0=\widetilde{\lambda}_0(\varepsilon) > \lambda_0$ in such a way that $\theta^2>2\widetilde{\lambda}_0$ in ${C^{\rm ext}_\varepsilon}$. The integral above converges and defines a continuous function of $\theta$ when $\theta^2>{2\widetilde{\lambda}_0}$. 
The last factor in the integrand can be estimated, for $\theta \gg 1$,  as follows
\[
e^{\widetilde{\lambda}_0\frac{|u|^2+4R\theta}{\theta^2}} -1 \leq \widetilde{\lambda}_0\frac{|u|^2+4R\theta}{\theta^2}
e^{\widetilde{\lambda}_0\frac{|u|^2+4R\theta}{\theta^2}} \leq C \frac{u^2 + 1}{\theta} e^{\frac{u^2}{4}}.
\]
Therefore, the integral behaves as $O(\theta^{-1})$ as $\theta\to \infty$, and thus
\begin{equation}\label{er1}
|r_1|\leq \frac {C(\varepsilon)}{|x-y|\theta}, \quad (t,x)\in {C^{\rm ext}_\varepsilon},~~t^2 +| x|^2\to\infty, ~~|y|\leq R.
\end{equation}

Similar arguments can be used to estimate $r_2$. Using (\ref{estp}) and the relation (\ref{rr2}) followed by substitution $z=y+\theta s\alpha+\sqrt su$, we  obtain
\[
|r_2|\leq C\int_{\mathbb{R}^d}\int_t^\infty \frac{1}{s^{d/2}}e^{-\frac{|\theta s\alpha+y-z|^2}{2s}+\widetilde{\lambda}_0s}|v(z)|ds dz
\]
\[=C\int_t^\infty\int_{\mathbb{R}^d}
e^{-\frac{|u|^2}{2}+\widetilde{\lambda}_0s}|v(y+\theta s\alpha+\sqrt su)| du ds.
\]
From (\ref{4}) it follows that $v$ in the integrand above vanishes when
\[
|u|^2<\sigma:=\theta^2s- 4R\theta. 
\]
Then
\[
|r_2|\leq C\int_t^\infty \int_{|u|^2>\sigma}e^{-\frac{|u|^2}{2}+\widetilde{\lambda}_0s} duds \leq C\int_t^\infty \eta
(\sigma) e^{-\frac{\theta^2s- 4R\theta}{2}+\widetilde{\lambda}_0s}
ds,
\]
where $\eta(\sigma)=\sigma^\frac{d-2}{2}$ if $\sigma\geq 1$ and $\eta(\sigma)=1$ if $\sigma< 1$.
We choose $\widetilde{\lambda}_0$ to be so close to $\lambda_0$ that $\theta^2/2-\widetilde{\lambda}_0>\delta$ with some $\delta>0$ when $(t,x)\in {C^{\rm ext}_\varepsilon}$. Since  the function $f(\theta):=\theta^2/2-\widetilde{\lambda}_0,~ \theta\geq \sqrt{2\lambda_0}+\varepsilon, $ is separated from zero and $f\sim \theta^2/2$ at infinity, there exists a constant $\delta_1>0$ such that $\theta^2/2-\widetilde{\lambda}_0>\delta_1\theta^2$. Since $\theta^2 t \gg R \theta$ for $(t,x)\in C^{\rm ext}_\varepsilon, ~t^2+|x|^2\to\infty, |y| \leq R$,
 the integral above can be estimated by $Ce^{-\delta_1\theta^2t/2}$. This provides the estimate on $r_2$, which, together with (\ref{er1}), implies the estimate of the remainder term in (\ref{tb}). In order to complete the proof of Part (b), it remains to justify the properties of function $a$: continuity, positivity, and the asymptotic behavior as $\theta\to\infty$.

Estimate (\ref{estp})   implies that the integral in (\ref{taaa}), with the integrand replaced by
its absolute value, can be estimated by
\[
 C \int_0^\infty \int_{\mathbb{R}^d} e^{-\frac{\theta^2s}{2}+ 2 R \theta+\tilde{\lambda}_0 s}|v(z)|
p_0(s,z- y)  d z ds,
\]
%where $\lambda_0 < \widetilde{\lambda}_0 <\theta^2/2-\delta(\sqrt{2 \lambda_0} + \varepsilon)^2 (1/2 - \delta_1)$
with $ \widetilde{\lambda}_0 =\lambda_0+\varepsilon^2/2 $, where $\varepsilon$ is the index in
$C^{\rm ext}_\varepsilon$. The latter integral
converges uniformly with respect to $\theta,\alpha, y$ when $\theta \geq \sqrt{2\lambda_0}+\varepsilon$. Hence $a$ is a continuous function of all its arguments.
%Similarly, the integral (\ref{taaa})
%converges uniformly after differentiation in $\theta,\alpha, y$, and therefore $a$ is infinitely smooth.

The existence of a lower positive bound for function $a$ is an obvious consequence of the following estimate justified  next:
\begin{equation} \label{apos}
p(t,x,y)\geq cp_0(t,x-y), ~~t\geq0, ~~|y|\leq R, ~~|x|\geq R, ~~c>0.
\end{equation}
If $|x|=R$ and $0\leq t\leq1$, then the validity of (\ref{apos}) with some $c=c_1$ follows from (\ref{apos1}) . Its validity with some $c=c_2$, when $|x|=R$ and $t\geq1$, follows from the facts that both functions $p$ and $p_0$ are continuous and positive there, $p\to\infty$ as $t\to\infty$, and $p_0$ vanishes as  $t\to\infty$. Hence,  (\ref{apos}) holds with  $c=\min(c_1,~c_2)$ when $|x|=R$. This implies  (\ref{apos}) for all $|x|\geq R$ due to the maximum principle for the heat equation in the region  $|x|\geq R$.

To justify the asymptotics of $a$ as $\theta\to\infty$, we split $a-1$ as $a-1=I_1+I_2$ where $I_1$ and $I_2$ are given by the same double integral (\ref{taaa}) with integration in $s$ restricted to $(0,\gamma)$ and $(\gamma,\infty)$, respectively, with $\gamma=(8R+1)/\theta$. One can check that the  function ${-\frac{\theta^2s}{2}+ 2 R \theta+\tilde{\lambda}_0 s}$  in the exponent of the integrand above does not exceed $-\theta^2s/4$ when $s>\gamma$ and $\theta$ is large enough. Hence
\[
|I_2|\leq C \int_\gamma^\infty \int_{\mathbb{R}^d} e^{-\theta^2s/4}|v(z)|
p_0(s,z- y)  d z ds\leq C_1 \int_\gamma^\infty e^{-\theta^2s/4} ds
\]
\begin{equation}\label{I22}
=4C_1\theta^{-2}e^{-(8R+1)\theta/4}, \quad~~~ \theta\to\infty,
\end{equation}
i.e., the contribution from $I_2$ to the asymptotics of $a-1$, as $\theta\to \infty$, can be disregarded.

We will use a better approximation of $p$ when studying $I_1$. Since $p\leq \hat{p}$, where $\hat{p}$ is the fundamental solution of problem (\ref{mpeq}) with $v(x)$ replaced by $\hat{v}:=\max v(x)$, and $\hat{p}(t,x,y)=e^{\hat{v}t}p_0(t,x -y)$, it follows that
\begin{equation} \label{apos1}
p(s,z,y)=p_0(s,z -y)(1+O(s)), \quad {\rm when} ~~~0\leq s\leq 1, ~~|y|\leq R, ~~z\in \mathbb{R}^d.
\end{equation}
We put it into the integral defining $I_1$ and then use (\ref{rr2}) and the substitution
 $z=y+\theta s\alpha+\sqrt su$. This leads to
 \[
 I_1=(2\pi)^{-d/2}\int_{\mathbb{R}^d}\int_0^\gamma
 e^{-\frac{|u|^2}{2}}v(y+\theta s\alpha+\sqrt su)(1+O(s))
 ds du.
 \]
 After substitution $s\to s/\theta$, we obtain
 \[
I_1=(2\pi)^{-d/2}\frac{1}{\theta}\int_{\mathbb{R}^d}\int_0^{8R+1}
 e^{-\frac{|u|^2}{2}}v(y+ s\alpha+u\sqrt {s/\theta})
 ds du+O(\theta^{-2}).
 \]

 If $v\in C^2$, then we write the integral above as a sum of three terms by splitting $v$ in the integrand as $v(y+ s\alpha)+\langle \nabla v (y+ s\alpha),u\sqrt {s/\theta}\rangle+O(|u|^2\theta^{-1})$. The first term coincides with the main term of asymptotics of $a$ since $v(y+ s\alpha)=0$ when $s>8R+1$. The second term vanishes since the integrand is odd in $u$. The last one can be combined with the remainder term in the asymptotics of $a$. The latter arguments when $v\in C^1$ or $v\in C$ are similar.  The asymptotic formula for $I_1$ together with (\ref{I22}) justify the asymptotic behavior (\ref{I22}) of $a$ as $\theta\to\infty$.

This completes the proof of Part (b).
\\

\noindent
{\it Proof of Part (c).}    First, we will prove the statement in Part (c)
when $(t,x)\in C_{\varkappa,\varepsilon_0}$ with arbitrary $\varepsilon_0>0$, where
\[
C_{\varkappa,\varepsilon_0} = \{(t,x):  \sqrt{2 (\lambda_0-\varkappa)} +\varepsilon_0 \leq \theta \leq \frac{1}{\varepsilon_0} \},
\]
and $\varkappa$ is the distance between $\lambda_0$ and the second largest eigenvalue of the operator $L  = \frac{1}{2} \Delta + v$, or $\varkappa=\lambda_0$ if operator $L  = \frac{1}{2} \Delta + v$ has only one eigenvalue $\lambda=\lambda_0$. Thus, for sufficiently small
$\varepsilon _0> 0$, the region
$C_{\varkappa,\varepsilon_0} $ contains a part of the interior region  $C^{\rm int}_0$, the conical surface $C$, and the exterior
region $C^{\rm ext}_0$  without a small conical neighborhood of the hyperplane $t=0$.

Denote by $\chi=\chi(t)$
the indicator function of the interval $(0,1)$, i.e $\chi(t)=1$ when $t\in(0,1)$, $\chi(t)=0$ when $t\notin(0,1)$. We rewrite (\ref{betaeq}) in the form
\begin{equation} \label{betaeqa}
p(t,x,y) =e^{\lambda_0 t}\psi(x)\psi(y)+ \beta_1(t,x,y)+ \beta_2(t,x,y), \quad \beta_1=\chi\beta, ~~ \beta_2=(1-\chi)\beta.
\end{equation}
Then we
split integral $I$ (see (\ref{duh}), (\ref{II})) as
\begin{equation} \label{tnti}
I=A+B_1+B_2
\end{equation}
 by substituting the sum (\ref{betaeqa}) for $p$ in (\ref{II}), and study each term separately.

Due to (\ref{betaeq1}),  for $\widetilde{\lambda}_0 > {\lambda}_0$, the following analogue of estimate (\ref{estp}) is valid for $\beta_2$:
\begin{equation} \label{estpa}
 |v(z)|\beta_2(s,z,y) \leq C(\widetilde{\lambda}_0) |v(z)|p_0(s,z-y)e^{\widetilde{(\lambda}_0 -\varkappa )s}, ~~s \geq 0,~z 
\in \mathbb{R}^d, ~|y| \leq R. 
\end{equation}

 The proof of Part (b) was based on estimate (\ref{estp}) on $p$. One can repeat the same arguments for integral $B_2
 $ using estimate (\ref{estpa}) on $\beta_2$ instead of the estimate (\ref{estp}) on $p$ (and using $\lambda_0-\varkappa$ instead of $\lambda_0$). This leads to the following statement. For each $\varepsilon_0 > 0$ and $(t,x) \in C_{\varkappa,\varepsilon_0}$,
\begin{equation} \label{tbb20}
B_2 =  p_0(t,x-y)
(a_{\beta_2}(\theta, \alpha, y)  + O(|x-y|^{-1})~~~{\rm as}~~|x| \rightarrow \infty,
\end{equation}
where $a_{\beta_2}$ is a continuous function equal to
\[
a_{\beta_2}(\theta, \alpha, y) =  \int_0^\infty \int_{\mathbb{R}^d} e^{-\frac{\theta^2s}{2}-\theta\langle \alpha, y-z \rangle}v(z)
\beta_2(s,z, y)  d z ds.
\]

Similar relations are valid for $B_1$. Indeed, by repeating arguments used to derive (\ref{ib}), we obtain
\begin{equation}\label{ib1}
B_1
=p_0(t,x-y)\int_{\mathbb{R}^d}\int_0^1(\frac{t}{ t-s})^{\frac{d}{2}} e^{-\frac{\theta^2s}{2}-\theta\langle \alpha, y-z \rangle-\frac{|\theta s\alpha+y-z|^2}{2(t-s)}} v(z) \beta_1(s,z,y)ds dz.
\end{equation}

 If $(t,x)\in {C}_{\varkappa,\varepsilon_0}, ~ |y|,|z|\leq R,$ and $|x|\to\infty$, then,  uniformly in $y,z$,
\[
\frac{t}{ t-s}=1+O(|x-y|^{-1}),~~~e^{-\frac{|\theta s\alpha+y-z|^2}{2(t-s)}}=1+O(|x-y|^{-1}).
\]
These relations, (\ref{ib1}) and integrability of $v(z)\beta_1$ imply an analogue of (\ref{tbb20}) for $B_1$. Hence for
$(t,x)\in {C}_{\varkappa,\varepsilon_0}$, we have
\begin{equation} \label{tbb}
B:=B_1+B_2 =  p_0(t,x-y)
(a_{\beta}(\theta, \alpha, y)  + O(|x-y|^{-1})~~~{\rm as}~~|x| \rightarrow \infty,
\end{equation}
where $a_{\beta}$ is a continuous function equal to
\begin{equation} \label{abeta}
a_{\beta}(\theta, \alpha, y) = \int_0^\infty \int_{\mathbb{R}^d} e^{-\frac{\theta^2s}{2}-\theta\langle \alpha, y-z \rangle}v(z)
\beta(s,z, y)  d z ds.
\end{equation}

It remains to study the asymptotic behavior at infinity of the integral
\begin{equation} \label{A}
A=\frac{\psi(y)}{(2 \pi)^{d/2}}  \int_{\mathbb{R}^d}\int_0^t\frac{1}{ (t-s)^{d/2}} e^{-\frac{|x-z|^2}{2(t-s)}+\lambda_0s} dsv(z) \psi(z) dz,
~~|y|\leq R.
\end{equation}

We will need the following lemma.
\begin{lemma}\label{xxx}
Let
\[
H(\omega)=\int_{-\infty}^l h(\tau)e^{-\omega\tau^2}d \tau, \quad   l\in \mathbb{R},
\]
where $h$ is $C^2$-smooth, $|h|<C e^{\tau^2}$, and the following relations are valid when $\tau\to-\infty$:
\[
 h\geq C|\tau|^{-2}, ~~|h'|\leq C|\tau| h.
\]
Then
\begin{equation}\label{Fl}
\begin{split}
H  & =\frac{\sqrt\pi}{2\sqrt \omega}(1+{\rm erf}(l\sqrt \omega))(h(0)+O(\omega^{-1})) \\
& +\frac{ h(0)-h(l)}{2 \omega l}e^{-\omega l^2}(1+O(\omega^{-1})),  \quad \omega\to\infty,
\end{split}
\end{equation}
where the estimates of the remainder terms are uniform in $l\in \mathbb{R}$.
\end{lemma}
\proof
Consider first the case when $l\geq 1$. Then the right-hand side in (\ref{Fl}) equals $\sqrt\frac{\pi}{\omega} (h(0)+O(\omega^{-1}))$, and the validity of (\ref{Fl}) is a simple consequence of the Laplace method. Let $l<1$. Since $h>0$ for $\tau\to-\infty$, and the statement is obviously valid when $h$ is a constant, it is enough to prove the statement when $l<1$ and $h$ is positive.

We represent $H$ in the form $H=H_0+\widetilde{H}$, where
\[
H_0=\int_{-\infty}^l h(0)e^{-\omega\tau^2}d\tau, ~ ~~\widetilde{H}=\int_{-\infty}^l [h(\tau)-h(0)]e^{-\omega\tau^2}d\tau.
\]
We divide and multiply the integrand in $\widetilde{H}$ by $-2\omega\tau$  and integrate by parts. This leads to
\[
\widetilde{H}= \frac{[h(0)-h(l)]}{2\omega l}e^{-\omega l^2}-\int_{-\infty}^l \frac{d}{d\tau}(\frac{[h(0)-h(\tau )]}{2\omega\tau})e^{-\omega\tau^2}d\tau, ~~\omega>1.
\]
From the properties of function $h$, it follows that the absolute value of the pre-exponential factor in the integrand above does not exceed $C\omega^{-1}h(\tau)$ (since $h$ is now considered to be positive), and therefore the last integral does not exceed $C\omega^{-1}H$. Hence
\[
H= H_0+\frac{[h( 0)-h(l)]}{2\omega l}e^{-\omega l^2}+O(\omega^{-1}H),  \quad \omega\to\infty.
\]
It remains to note that $H_0$ coincides with the first term in the right-hand side of (\ref{Fl}).
\qed
\\

 Denote by $F$
the interior integral in (\ref{A}). We will find its asymptotics as $|x|\to\infty,~|z|\leq R,$ considering $F$ as a function of $\theta'=\frac {|x-z|}{t},~|z|\leq R$, and then we will express $\theta'$ through $\theta$.
 After the substitution $t-s=|x-z|\sigma$, integral $F$ takes the form
\[
F=|x-z|^{1 - \frac{d}{2}}e^{\lambda_0 t}\int_0^{1/\theta'}\frac{1}{\sigma^{d/2}}e^{- |x-z|(\frac{1}{2\sigma}+ \lambda_0\sigma )}d\sigma, ~~|z|\leq R.
\]
Here the function $f(\sigma)=\frac{1}{2\sigma}+ \lambda_0\sigma $, considered on the semi-axis $\sigma>0$, has a single stationary point (minimum) at $\sigma=1/\sqrt{2\lambda_0}$, and the asymptotic behavior of the integral $F$ as $|x-z|\to\infty$ depends essentially on whether the limit of integration $1/\theta'$ is greater or smaller than $1/\sqrt{2\lambda_0}$, and how close it is to that value.

We take $e^{-\sqrt{2\lambda_0}|x-z|}$ (the value of the exponential function at the stationary point) out of the integral, and make the substitution $\tau=\frac{\sqrt{2\lambda_0}\sigma-1}{\sqrt{2\sigma}}$, which makes the exponent of the integrand equal to $-|x-z|\tau^2$:

\begin{equation}\label{Ig}
F=|x-z|^{1-\frac{d}{2}}e^{\lambda_0 t-\sqrt{2\lambda_0}|x-z|}\int_{-\infty}^{g(\theta')} h(\tau)e^{-|x-z|\tau^2}d\tau, \quad g(\theta')=\frac{\sqrt{2\lambda_0}-\theta'}{\sqrt{2\theta'}},
\end{equation}
where $|z|\leq R,~h(\tau)=\frac{ \sigma'(\tau)}{[\sigma(\tau)]^{d/2}}\in C^\infty,$
\begin{equation}\label{h0}
h(0)=\sqrt{2}(2\lambda_0)^\frac{d-3}{4},
\end{equation}
 and $h$ behaves as a power function at infinity (one can show that $h(\tau)\sim c_-|\tau|^{d-3}$ as $\tau\to-\infty$, $h(\tau)\sim c_+\tau^{1-d}$ as $\tau\to\infty$, where $c_\pm>0$).

We apply Lemma \ref{xxx} with $l=g(\theta')$ and $\omega =|x-z|$ to (\ref{Ig}), and obtain
\[
F=\frac{\sqrt\pi}{2}\omega^{\frac{1- d}{2}}(1+{\rm erf}(\sqrt{t}\frac{\sqrt { 2\lambda_0}-\theta'}{\sqrt { 2 }}))e^{\lambda_0 t-\sqrt{2\lambda_0}\omega}(h(0)+O(\omega^{-1}))
\]
\[
+\omega^{-\frac{d}{2}}\frac{h(0)-h(g(\theta'))}{2g(\theta')}e^{-\frac{\omega^2}{2t}}(1+O(\omega^{-1})),~~|z|\leq R.
\]
We put $\omega=|x|-\langle \dot{x}, z\rangle+O(|x-y|^{-1})$ in the first term on the right and $\omega=|x-y|+\langle \alpha, y-z\rangle+O(|x-y|^{-1})$ in the second term. This leads to
\[
F=\frac{\sqrt\pi}{2}|x|^{\frac{1- d}{2}}(1+{\rm erf}(\sqrt{t}\frac{\sqrt { 2\lambda_0}-\theta'}{\sqrt { 2 }}))e^{\lambda_0 t-\sqrt{2\lambda_0}(|x|-\langle\dot x,z\rangle)}(h(0)+O(|x|^{-1/2}))
\]
\begin{equation}\label{C}
+(2\pi)^{d/2}p_0 (t,x-y)
\theta^{-\frac{d}{2}} \frac{h(0)-h(g(\theta'))}{2g(\theta')} e^{-\theta\langle \alpha, y-z \rangle}(1+O(|x-y|^{-1})),~~|y|,|z|\leq R.
\end{equation}
 We replace everywhere $\theta'$ by $\theta+O(|x-y|^{-1})$. Note that the ratio  in the second term is smooth in $\theta'$ since the denominator vanishes only at $\theta'=\sqrt{2\lambda_0}$ and has zero of the first order there (see (\ref{Ig})), and the numerator also vanishes at that point. Thus
\begin{equation}\label{C1}
\frac{h(0)-h(g(\theta'))}{2g(\theta')}=\frac{h(0)-h(g(\theta))}{2g(\theta)}+O(|x-y|^{-1}).
\end{equation}
Let
 \[
 u':=\sqrt{t}\frac{\sqrt { 2\lambda_0}-\theta'}{\sqrt { 2 }},~~u:=\sqrt{t}\frac{\sqrt { 2\lambda_0}-\theta}{\sqrt { 2 }},~~\Delta u=u'-u.
 \]
 %Then $|\Delta u|=O(\frac{1}{\sqrt{t}})=O(\frac{1}{\sqrt{|x|}})$ when $(t,x) \in C_{\varkappa, \varepsilon_0},~ |x|\to\infty$. Therefore,

If
 $u\geq -1$ and $|\Delta u|\leq 1/2$, then
\begin{equation}\label{CCC}
 2\geq(1+{\rm erf}(\sqrt{t}\frac{\sqrt { 2\lambda_0}-\theta}{\sqrt { 2 }})\geq c_0>0,
 \end{equation}
  and
\begin{equation}\label{C2}
 1+{\rm erf}(\sqrt{t}\frac{\sqrt { 2\lambda_0}-\theta'}{\sqrt { 2 }})=(1+{\rm erf}(\sqrt{t}\frac{\sqrt { 2\lambda_0}-\theta}{\sqrt { 2 }}))(1+O(\frac{1}{\sqrt{|x|}})),~~u \geq -1.
\end{equation}
A different formula is valid in the case when $u\leq-1$. Observe that
\[
1+ {\rm erf}(u)=|\sqrt{\pi}u|^{-1}e^{-u^2}g(u),~~u \leq -1,
\]
where $g(u)=1+O(|u|^{-1})$ as  $ u \rightarrow -\infty$, $g$ is separated from zero, and its  derivative is bounded. Hence,
\[
g(u')=g(u)+O(\Delta u)=g(u)(1+O(\Delta u)) ~~ {\rm when} ~~u\leq -1,~|\Delta u|\leq 1/2.
\]
Obviously, when $u\leq -1,~|\Delta u|\leq 1/2$, we have
\[
|u'|^{-1}=|u|^{-1}(1+O(\Delta u)) ,  ~~ e^{-u'^2}=e^{-u^2-2u\Delta u}(1+O(\Delta u)),
\]
and therefore,
\[
1+ {\rm erf}(u')=(1+ {\rm erf}(u))e^{-2u\Delta u}(1+O(\Delta u)), ~~u\leq -1, ~~|\Delta u|\leq 1/2.
\]
Since $\Delta u=\sqrt\frac{1}{2t}\langle\dot{x},z-y\rangle+O(\frac{1}{|x|\sqrt t})$ when
 $u\leq-1,~(t,x) \in C_{\varkappa, \varepsilon_0}, ~|x|\to\infty$, it follows that
\begin{equation}\label{C3}
(1+{\rm erf}(\sqrt{t}\frac{\sqrt { 2\lambda_0}-\theta'}{\sqrt { 2 }}))=(1+{\rm erf}(\sqrt{t}\frac{\sqrt { 2\lambda_0}-\theta}{\sqrt { 2 }}))e^{(\theta-\sqrt { 2\lambda_0})\langle \dot{x},z-y \rangle}(1+O(\frac{1}{\sqrt{|x|}})).
\end{equation}
Due to (\ref{CCC}), formulas (\ref{C2}), (\ref{C3}) are equivalent when $u\in [-1,0]$ since the latter inclusion implies that $0\leq\theta-\sqrt { 2\lambda_0}\leq \sqrt{2/ t}=O(1/\sqrt {|x|}),~|x|\to\infty $. Hence, from (\ref{C2}) and (\ref{C3}), it follows that
\begin{equation}\label{C5}
(1+{\rm erf}(\sqrt{t}\frac{\sqrt { 2\lambda_0}-\theta'}{\sqrt { 2 }}))=
(1+{\rm erf}(\sqrt{t}\frac{\sqrt { 2\lambda_0}-\theta}{\sqrt { 2 }}))b(\theta,\dot{x},y,z)(1+O(\frac{1}{\sqrt{|x|}}))
\end{equation}
when $(t,x) \in C_{\varkappa, \varepsilon_0},~ |x|\to\infty$, where $b=1 $ if $\theta\leq \sqrt { 2\lambda_0},~
b=e^{(\theta-\sqrt { 2\lambda_0})\langle\dot{x},z-y\rangle}$ if $\theta\geq \sqrt { 2\lambda_0}$.

We put (\ref{C1}),(\ref{C5}) into (\ref{C}) and combine the resulting formula for $F$ with (\ref{A}). This implies
\[
A=\psi(y)|x|^{\frac{1- d}{2}}e^{\lambda_0 t-\sqrt{2\lambda_0}|x|}(1+{\rm erf}(\sqrt{t}\frac{\sqrt { 2\lambda_0}-\theta}{\sqrt { 2 }}))(\gamma_1(\theta,\dot{x},y)+O(|x|^{-1/2}))
\]
\begin{equation}\label{AC}
+p_0 (t,x-y)
(\gamma_2(\theta,\alpha,y)+O(|x-y|^{-1})):=A_1+A_2,~~|y|\leq R,
\end{equation}
where
\begin{equation} \label{ga1}
\gamma_1=(2\pi)^\frac{1-d}{2}2^{-3/2}h(0)\int_{|z|\leq R}e^{\sqrt{2\lambda_0}\langle\dot x,z\rangle}b(\theta,\dot{x},y,z)v(z)\psi(z)dz,
\end{equation}
\begin{equation} \label{ga2}
\gamma_2=
\frac{\psi(y)}{(2 \pi)^{d/2}} \theta^{-\frac{d}{2}} \frac{h(0)-h(g(\theta))}{2g(\theta)}\int_{|z|\leq R} e^{-\theta\langle \alpha, y-z \rangle}v(z)\psi(z)dz.
\end{equation}
Using (\ref{h0}), one can easily check that $\gamma_1=\frac{1}{2}C(\dot x)$ when $\theta\leq\sqrt{2\lambda_0}$ where $C(\dot x)$ is defined in (\ref{ppsi}), (\ref{intrep}). Hence, for $(t,x)\in C_{\varkappa,\varepsilon_0}$, the first term in (\ref{AC}) can be rewritten in the form
\begin{equation}\label{A1C}
A_1=e^{\lambda_0t}\psi(x)\psi(y)(1+{\rm erf}(\sqrt{t}\frac{\sqrt { 2\lambda_0}-\theta}{\sqrt { 2 }}))(a_1(\theta,\dot x,y)+O(|x|^{-1/2})),
\end{equation}
where $a_1=\frac{1}{2}$ when $\theta\leq\sqrt{2\lambda_0}, ~ a_1=[C(\dot x)]^{-1}\gamma_1$ when $\theta\geq\sqrt{2\lambda_0}$.

We combine formulas (\ref{AC}), (\ref{A1C}) for $A$ with (\ref{duh}), (\ref{tnti}), and (\ref{tbb}), and obtain (\ref{tc})
with the coefficient $a_1$ defined in (\ref{A1C}) and $a_2=1+\gamma_2+a_\beta$.
This completes the proof of  (\ref{tc}) when $\sqrt{2(\lambda_0-\varkappa)}+\varepsilon_0\leq\theta\leq 1/\varepsilon_0$.

We extend formula (\ref{tc}) for  $0\leq\theta\leq\sqrt{2(\lambda_0-\varkappa)}+\varepsilon_0$ using an arbitrary continuous extension of $a_2$ to that region and using $a_1=1/2$ there.  Let us show that this extension provides a correct formula for $p$. In fact, we will justify a stronger statement that the main term of asymptotics in (\ref{tc}) coincides with the one in (\ref{tt}) when $ (t,x)\in C_{\varepsilon}^{\rm int}$. Indeed, the function erf in (\ref{tc}) is exponentially, in $t$, close to one, and the first term of the asymptotics in (\ref{tc}) coincides with the one in (\ref{tt}). The second term in (\ref{tc}) is exponentially smaller than the first one and can be omitted when $ (t,x)\in C_{\varepsilon}^{\rm int}$. This can be easily checked by comparing  $p(t,x-y)$ and the exponential factors in (\ref{tt}), (\ref{ppsi}).
\qed
\section{Appendix} \label{Att}

Here we discuss  the relation between main terms of asymptotics in (\ref{tc}) in a neighborhood
of the cone $C$ and  the relation between (\ref{tb}) and (\ref{tc}) in $C_{\varepsilon}^{\rm ext}$.

 Since
\[
e^{\lambda_0t-\sqrt{2\lambda_0}|x|} =e^{\lambda_0t-\sqrt{2\lambda_0}|x-y|}(e^{-\sqrt{2\lambda_0}<\dot x,y>}+O(\frac{1}{|x-y|}))= e^{-
\frac{|x-y|^2}{2t}}(e^{-\sqrt{2\lambda_0}<\dot x,y>}+O(\frac{1}{|x-y|}))
\]
when $\theta =\sqrt{2\lambda_0}$,
it follows from (\ref{ppsi}) that the exponentially growing factors in the main terms of asymptotics in (\ref{tc})
coincide on $C$. The power factors are $|x|^{\frac{1-d}{2}}$ in the first term and $t^{\frac{-d}{2}}$ in the second term.
Taking into account that the erf-function is zero on $C$ and $a_1 = 1/2$ there, we obtain
that the second term on $C$ is smaller by a factor of order $O(t^{-1/2}), ~t \to\infty$. This statement
remains valid in a neighborhood of $C$ where the argument of the function erf is bounded.

Consider the region $\sqrt t(\theta-\sqrt{2\lambda_0})\geq1$ (where the argument of the erf function in (\ref{tc})
is separated from zero). In this region, the first term in (\ref{tc}) can be rewritten in the form
$d(\theta,\alpha, y)p_0(t, x - y)(1 + O(|x - y|^{-1/2}))$ with a certain continuous function $d$. The latter
formula with an explicit form of $d$ can be obtained from (\ref{ppsi}), the expression for $a_1$ (after formula (\ref{A1C})), and the asymptotics of the function erf:
\begin{equation}\label{aerf}
1 + {\rm erf}(u) = |\sqrt \pi u|^{-1}e^{-u^2}
(1 + O(|u|^{-1})),~ u \to\ -\infty.
\end{equation}
In particular, $c_1/(\theta-\sqrt{2\lambda_0}) \leq d(\theta,\alpha, y) \leq c_2/(\theta-\sqrt{2\lambda_0}) $ with some $c_1, c_2 > 0$  for
sufficiently small $\theta-\sqrt{2\lambda_0}$. Therefore, the first term in (\ref{tc}) is much larger than the second one if $\theta  - \sqrt{2 \lambda_0}$ is positive and small.

 Now let $(t,x) \in C_{\varepsilon}^{\rm ext}$ and $\theta \leq 1/\varepsilon_0$.
 It was mentioned at the end of the proof of the theorem that the coefficient $a_2$ in (\ref{tc}) is equal to
 $1 + \gamma_2+a_\beta $, where $\gamma_2$ is defined in
(\ref{ga2}) and  $a_\beta$ is given in (\ref{abeta}). Using the expressions for $h$ and $g$ from (\ref{Ig}) and (\ref{h0}), one can specify $\gamma_2$ as follows:
\[
\gamma_2=[\frac{-(\theta/\sqrt{2\lambda_0})^{\frac{1-d}{2}}}{\sqrt{2\lambda_0}(\sqrt{2\lambda_0}-\theta)}+
\frac{2}{ {\theta^2} -2 \lambda_0}]\psi(y)
\int_{\mathbb{R}^d} e^{-\theta\langle \alpha, y-z \rangle} v(z) \psi(z) dz.
\]
The first factor in the right-hand  side of this formula   is a sum of two terms.
Each of them has a singularity at $\theta=\sqrt{2\lambda_0}$, but the sum is smooth.

%If $(t,x) \in C_{\varepsilon}^{\rm ext}$ (i.e., $\theta >\sqrt {2 \lambda_0}+\varepsilon$), then the first term in
%(\ref{tc}) can be rewritten in the form $d(\theta,\alpha,y)p_0(t,x-y)(1+O(|x-y|^{-1/2}))$ with a certain continuous function $d$. The latter %formula with an explicit form of $d$ can be obtained from (\ref{ppsi}) and the asymptotics of the function erf:
%\begin{equation}\label{aerf}
%1+ {\rm erf}(u)=|\sqrt{\pi}u|^{-1}e^{-u^2}(1+O(|u|^{-1})),~~u\to-\infty.
%\end{equation}
%Function
An explicit calculation of the coefficient $d$ shows that it is given by the same expression as $-\gamma_2$ with the second term of the first factor in $\gamma_2$ omitted. Hence, (\ref{tc}) in $C_{\varepsilon}^{\rm ext}$ takes the form
\[
p(t,x,y) = p_0(t,x-y)(1+a_\beta + \frac{2}{ \theta^2-2 \lambda_0}\psi(y)
\int_{\mathbb{R}^d} e^{-\theta\langle \alpha, y-z \rangle} v(z) \psi(z) dz+O(|x-y|^{-\frac{1}{2}})).
\]
If we replace here $ 2/({\theta^2}-2\lambda_0)$ by $\int_0^\infty e^{-\theta^2s/2+\lambda_0s}ds$ and substitute  (\ref{abeta}) for $a_\beta$, this formula for $p$ will coincide with formula  (\ref{tb}), where $a$ is given by (\ref{taaa}). This provides a direct justification of the equality of the main terms of asymptotics in Parts (b) and (c) when $\theta \geq\sqrt {2 \lambda_0}+\varepsilon$.

The main reason to provide this justification was to show that the two main terms of asymmptotics of $p$ in formula (\ref{tc}) have the same order when $(t,x)\in C_{\varepsilon}^{\rm ext}$, since an indirect justification of the equivalency of main terms in (\ref{tb}) and (\ref{tc}) is obvious: both formulas give the asymptotics of the same function $p$.
\\
\\
\noindent {\bf \large Acknowledgments}:
The work of  L. Koralov was supported by the Simons Foundation Fellowship (award number 678928)
and by the Russian Science Foundation,  project ${\rm N}^o$~20-11-20119.
The work of B. Vainberg was supported by the Simons Foundation grant 527180.

\end{document}